\documentclass[11pt]{elsarticle}
\usepackage{amsfonts,amsmath,amssymb,dsfont}
\usepackage{amsfonts,amsmath}
\usepackage{hyperref}
\usepackage{graphicx}
\usepackage{subfigure}

\usepackage{tikz}
\usetikzlibrary{calc,arrows}

\newcommand{\detail}[1]{\par\noi{\bf [Proof detail\ }{#1}
\hfill{\bf ]}\par\noi\hspace{-4pt}}
\renewcommand{\detail}[1]{}

\newcommand{\dis}{\displaystyle}

\newcommand{\med}{\medskip}
\newcommand{\noi}{\noindent}

\newcommand{\halmos}{\rule{1ex}{1.4ex}}

\newcommand{\quand}{\quad\mbox{and}\quad}

\newtheorem{theorem}{Theorem}
\newtheorem{proposition}[theorem]{Proposition}
\newtheorem{corollary}[theorem]{Corollary}
\newtheorem{conjecture}[theorem]{Conjecture}
\newtheorem{lemma}[theorem]{Lemma}

\newtheorem{remark}[theorem]{Remark}

\newcommand{\bt}{\begin{theorem}}
\newcommand{\et}{\end{theorem}}
\newcommand{\bl}{\begin{lemma}}
\newcommand{\el}{\end{lemma}}
\newcommand{\bp}{\begin{proposition}}
\newcommand{\ep}{\end{proposition}}
\newcommand{\bcor}{\begin{corollary}}
\newcommand{\ecor}{\end{corollary}}
\newcommand{\br}{\begin{remark}\rm}
\newcommand{\er}{\end{remark}}
\newcommand{\bcon}{\begin{conjecture}}
\newcommand{\econ}{\end{conjecture}}

\newcommand{\QED}{\nopagebreak{\hspace*{\fill}$\halmos$\medskip}}
\newenvironment{Proof}[1][]{\noi\textbf{Proof #1}}{\QED}
\newcommand{\bpro}{\begin{Proof}}
\newcommand{\epro}{\end{Proof}}

\newcommand{\be}{\begin{equation}}
\newcommand{\ee}{\end{equation}}
\newcommand{\ba}{\begin{array}}
\newcommand{\ea}{\end{array}}
\newcommand{\bc}{\be\begin{array}{r@{\,}c@{\,}l}}
\newcommand{\ec}{\end{array}\ee}

\newcommand{\de}{\delta}

\newcommand{\eps}{\varepsilon}
\newcommand{\la}{\lambda}

\newcommand{\sig}{\sigma}

\newcommand{\Hi}{{\cal H}}

\newcommand{\R}{{\mathbb R}}
\newcommand{\N}{{\mathbb N}}
\newcommand{\Z}{{\mathbb Z}}

\newcommand{\E}{{\mathbb E}}
\renewcommand{\P}{{\mathbb P}}

\newcommand{\Yb}{{\mathbf Y}}

\newcommand{\volgt}{\ensuremath{\Rightarrow}}
\newcommand{\up}{\uparrow}
\newcommand{\down}{\downarrow}
\newcommand{\sub}{\subset}
\newcommand{\beh}{\backslash}

\newcommand{\asto}[1]{\underset{{#1}\to\infty}{\longrightarrow}}
\newcommand{\Asto}[1]{\underset{{#1}\to\infty}{\Longrightarrow}}

\newcommand{\Ato}[2]{\underset{{#1}\to{#2}}{\Longrightarrow}}

\newcommand{\ti}{\tilde}

\newcommand{\un}{\underline}

\newcommand{\ffrac}[2]{{\textstyle\frac{{#1}}{{#2}}}}
\newcommand{\dif}[1]{\ffrac{\partial}{\partial{#1}}}

\newcommand{\di}{\mathrm{d}}
\newcommand{\half}{{[0,\infty)}}
\newcommand{\expo}{\mbox{\large\it e}}
\newcommand{\ex}[1]{\expo^{\,\textstyle{#1}}}

\newcommand{\var}{{\rm Var}}

\newcommand{\ha}{\ffrac{1}{2}}

\begin{document}

\renewcommand{\labelenumi}{{\rm(\roman{enumi})}}

\title{The limiting shape of a full mailbox}

\author[mf]{M. Formentin
}
\ead{marco.formentin@rub.de}

\author[mf]{J.M. Swart\corref{cor}}
\ead{swart@utia.cas.cz}\cortext[cor]{Corresponding author}
\ead[url]{http://staff.utia.cas.cz/swart/}

\address[mf]{Institute of Information Theory and Automation of the ASCR, Pod Vod\'{a}renskou v\v{e}\v{z}\'{i} 4, CZ-18208 Prague, Czech Republic}

\begin{keyword}
Self-organized criticality, Gabrielli and Caldarelli queueing model, Barabasi
queueing model, email communication.
\MSC[2010] 82C27, 60K35, 82C26, 60K25
\end{keyword}
%

\date{\today}

\begin{abstract}
We study a model for email communication due to Gabrielli and Caldarelli,
where someone receives and answers emails at the times of independent Poisson
processes with intensities $\la_{\rm in}>\la_{\rm out}$. The receiver assigns
i.i.d.\ priorities to incoming emails according to some atomless law and
always answers the email in the mailbox with the highest priority. Since
the frequency of incoming emails is higher than the frequency of answering,
below a critical priority, the mailbox fills up ad infinitum. We prove a
theorem about the limiting shape of the mailbox just above the critical point,
linking it to the convex hull of Brownian motion. We conjecture that this
limiting shape is universal in a class of similar models, including a model
for the evolution of an order book due to Stigler and Luckock.
\end{abstract}

\maketitle

\section{Introduction}

\subsection{Description of the model}\label{S:descrip}

The queueing model we study in this paper was introduced by Gabrielli and
Caldarelli \cite{CG09} as a variation of Barab\'asi's queueing system
\cite{Bar05}. These and similar models have been investigated in recent years
in the complex system literature as they are able to capture some universal
patterns in human written communication \cite{Bla07,MF14,Vaz06}. Usually in
this context, the model aims to describe the response time statistics of a
user which assigns a subjective priority to each incoming message (email,
paper mail or sms) and then answers first the highest priority one. We will be
interested in a somewhat different function of the process, namely, the
asymptotics of the number of items waiting to be executed with priority close
to the critical point. The mathematical description of the model is as
follows.

Tasks arrive according to a Poisson point process with rate $\la_{\rm
  in}$. Each incoming task is assigned a priority. The priorities of incoming
tasks are i.i.d.\ real-valued random variables with some law $\mu$. We assume
that $\mu$ is atomless, which assures that all tasks in the queue have a
different priority. At times of a Poisson point process with rate $\la_{\rm
  out}$, the task with the highest priority in the queue is executed (i.e.,
removed from the queue). If at such a time, the queue is empty, then nothing
happens.

Since only the relative order of the priorities matters, the precise choice of
the law $\mu$ does not matter. For definiteness, we choose for $\mu$ the
uniform distribution on $[-\la_{\rm in},0]$. Note that our priorities are
negative numbers, i.e., $0$ is the highest possible priority, which will be
convenient from a mathematical perspective. By time scaling, we can without
loss of generality assume that $\la_{\rm out}=1$ so that our model depends on
a single parameter $\la:=\la_{\rm in}$.

Let $\Pi^\la_{\rm in}$ be the random collection of all pairs $(p,t)$, where
$t\geq 0$ is a time when a new task arrives and $p\in[-\la,0]$ is the priority
assigned to this task. Also, let $\Pi_{\rm out}$ be the collection of all
times $s\geq 0$ when tasks are executed. Let $Y^\la(t)$ be a finite subset of
$[-\la,0]$, describing the priorities of all tasks in the queue at time
$t$. Since all tasks have a different priority, the cardinality
$N^\la(t):=|Y^\la(t)|$ of $Y^\la(t)$ equals the number of tasks in the
queue. By convention, we choose $Y^\la(t)$ right-continuous in $t$
and let $Y^\la(t-):=\lim_{s\up t}Y^\la(s)$ denote the state of $Y^\la$
immediately prior to $t$. We order the elements of $Y^\la(t)$ from the highest
to the lowest priority:
\be\label{priord}
Y^\la(t)=\{Y^\la_1(t),\ldots,Y^\la_{N^\la(t)}\}
\quad\mbox{with}\quad
Y^\la_1(t)>Y^\la_2(t)>\cdots>Y^\la_{N^\la(t)}(t).
\ee
We start the process with $Y^\la(0)$ some finite subset of $[-\la,0]$.
Then $(Y^\la(t))_{t\geq 0}$ is a continuous-time Markov process
with the following description:
\begin{enumerate}
\item For each $(p,t)\in\Pi^\la_{\rm in}$, at time $t$, the previous state
  $Y^\la(t-)$ of the process is replaced by $Y^\la(t):=Y^\la(t-)\cup\{p\}$.
\item For each $t\in\Pi_{\rm out}$, at time $t$, the previous state
  $Y^\la(t-)$ of the process is replaced by
  $Y^\la(t):=Y^\la(t-)\beh\{Y^\la_1(t-)\}$ if $Y^\la(t-)\neq\emptyset$, and
  stays empty otherwise.
\end{enumerate}
We call $Y^\la=(Y^\la(t))_{t\geq 0}$ the \emph{inbox process} with
\emph{rate of incoming tasks} $\la\geq 0$.

We observe that $\Pi^\la_{\rm in}$ is a Poisson process on
$[-\la,0]\times\half$ with intensity one. It will be convenient to construct
$\Pi^\la_{\rm in}$ in the following way: letting $\Pi_{\rm in}$ denote a
Poisson process on $(-\infty,0]\times\half$ with intensity one, we define
$\Pi^\la_{\rm in}$ as the restriction
\be
\Pi^\la_{\rm in}:=\Pi_{\rm in}\cap\big([-\la,0]\times\half\big).
\ee
We observe the following consistency relation. If $0\leq\la'\leq\la$, then
setting
\be
Y^{\la'}(t):=Y^\la(t)\cap[-\la',0]\qquad(t\geq 0)
\ee
is exactly the inbox process with rate of incoming tasks $\la'$, started in
the initial state $Y^{\la'}(0):=Y^\la(0)\cap[-\la',0]$ and constructed from
the Poisson processes $\Pi^{\la'}_{\rm in}$ and $\Pi_{\rm out}$. Indeed,
incoming tasks with a priority below $-\la'$ have no influence on $Y^{\la'}$.
Also, at times when a task with priority below $-\la'$ is executed, the random
set $Y^{\la'}(t)$ is empty and stays empty, in line with the rules above.

In view of this, we can remove the last free parameter of our model and,
starting from a locally finite\footnote{By definition, a subset $Y$ of a
  topological space $X$ is \emph{locally finite} if $Y\cap C$ is finite for
  each compact subset $C$ of $X$.} subset $Y(0)\sub(-\infty,0]$, define an
``infinite'' process $(Y(t))_{t\geq 0}$ taking values in the locally finite
subsets of $(-\infty,0]$ such that for each $\la\geq 0$,
\be\label{Yla}
Y^\la(t)=Y(t)\cap[-\la,0]\qquad(t\geq 0)
\ee
is the inbox process with rate of incoming tasks $\la$, started in the initial
state $Y^\la(0):=Y(0)\cap[-\la,0]$. Formally, the process $Y$ follows the same
rules as $Y^\la$, with $\Pi^\la_{\rm in}$ replaced by $\Pi_{\rm in}$. Because
of consistency, $Y$ is well-defined, even though the set of times $\{t\geq
0:(p,t)\in\Pi_{\rm in}\}$ is a.s.\ dense in $\half$.

\subsection{The critical point}

Recall that $N^\la(t)$ denotes the number of tasks in the queue with priority
in $[-\la,0]$. We observe that $N^\la=(N^\la(t))_{t\geq 0}$ is a continuous-time
random walk with reflection at the origin, i.e., $N^\la$ is a Markov process
with state space $\N$ that jumps
\be\label{driftwalk}
n\mapsto n+1\mbox{ with rate }\la
\quand
n\mapsto n-1\mbox{ with rate }1_{\{n>0\}}.
\ee
This process is positive recurrent for $\la<1$, null recurrent for $\la=1$,
and transient for $\la>1$. We order the elements of $Y(t)$ as
  $Y_1(t)>Y_2(t)>\cdots$ as in (\ref{priord}). Transience for $\la>1$ and
recurrence for $\la\leq 1$ imply that
\be
\liminf_{t\to\infty}Y_1(t)=-1\quad{\rm a.s.},
\ee
so for each $\la>1$, there is a random time after which no tasks with a
priority below {$-\la$} are executed anymore. On the other hand, for each
$\la<1$, positive recurrence implies that the highest priority $Y_1(t)$ in the
inbox spends a positive fraction of time below $-\la$. In view of this, the
following proposition should not come as a surprise. Recall that a subset of a 
topological space is \emph{locally finite} if its intersection with any
compact set is finite. In particular, a subset $Y\sub(-1,0]$ is locally finite
if $Y\cap[-\la,1]$ is finite for all $0<\la<1$.

\bp[Long-time limit]\label{P:ergod}
There exists a random, locally finite subset
$Y^1(\infty)\sub(-1,0]$ such that, regardless of the initial state $Y(0)$,
\be\label{ergod}
\P\big[Y^\la(t)\in\cdot\,\big]
\asto{t}\P\big[Y^1(\infty)\cap[-\la,0]\in\cdot\,\big]
\qquad(0<\la<1),
\ee
where $\to$ denotes convergence of probability measures in total variation
norm distance. The random point set $Y^1(\infty)$ a.s.\ has infinitely
many elements. Writing
\be
Y^\la(\infty):=Y^1(\infty)\cap[-\la,0]
\quand
N^\la(\infty):=|Y^\la(\infty)|
\qquad(0\leq\la\leq 1),
\ee
one has
\be
\P\big[N^\la(\infty)=n\big]=\dis(1-\la)\la^n
\qquad(0\leq\la<1,\ n\geq 0).
\ee
\ep

The random set $Y^1(\infty)$ describes the long-time limit of the collection
of all tasks in the inbox with priorities above the critical point $-\la_{\rm
  c}:=-1$ waiting to be executed. We are interested in the shape of
$Y^1(\infty)$ near the critical point.

\subsection{The limiting shape near the critical point}\label{S:Hi}

Recall from Proposition~\ref{P:ergod} that $N^\la(\infty)$
denotes the equilibrium number of tasks with priority $>-\la$ in the inbox,
which is a.s.\ finite by the local finiteness of $Y^1(\infty)\sub(-1,0]$. We
will be interested in the shape of the random function $\la\mapsto
 N^\la(\infty)$ in the vicinity of $\la_{\rm c}=1$. To this aim, for
$\eps>0$, we define
\be\label{Heps}
H^\eps_s:=\eps N^{1-2\eps s}(\infty)\qquad(s>0).
\ee
Clearly, $H^\eps:(0,\infty)\to\half$ is a.s.\ right-continuous, nonincreasing,
$H^\eps(s)=0$ for $s\geq 1/(2\eps)$, and $\lim_{s\down 0}H^\eps(s)=\infty$ by the
fact (proved in Proposition~\ref{P:ergod}) that $|Y^1(\infty)|=\infty$ a.s.
As $\eps\down 0$, the function $H^\eps$ describes the shape of the function
$\la\mapsto N^\la(\infty)$ for $\la$ just below $\la_{\rm c}=1$, where we scale
distances in $\la$ by a factor $(2\eps)^{-1}$ and at the same time scale down
the numbers $N^\la(\infty)$ by giving each task a weight $\eps$.

We will prove that the random function $H^\eps$ converges as $\eps\down 0$ to
a random limiting function, that is closely linked to the convex hull of
Brownian motion. To formulate this properly, let $\Hi$ denote the space
of all right-continuous, nonincreasing functions $h:(0,\infty)\to\half$.
Then $h\in\Hi$ if and only if $h$ is the distribution function of a locally
finite measure on $(0,\infty]$, i.e., each $h\in\Hi$ corresponds to a locally
finite measure $\mu$ on $(0,\infty]$ such that
\be\label{hmu}
h_s=\mu\big((s,\infty]\big)\qquad\big(s\in(0,\infty)\big).
\ee
We equip $\Hi$ with the topology of vague convergence of the corresponding
locally finite measures on $(0,\infty]$. The following theorem is our main
result.

\bt[Limiting shape near the critical point]\label{T:main}
One has
\be\label{toH}
\P\big[(H^\eps_s)_{s>0}\in\,\cdot\,\big]
\underset{\eps\down 0}{\Longrightarrow}
\P\big[(H_s)_{s>0}\in\,\cdot\,\big],
\ee
where $\Rightarrow$ denotes weak convergence of probability laws on $\Hi$ with
respect to the topology defined above, and $H_s:=\sup_{t\geq
  0}\left(B_t-st\right)$ with $\left(B_t\right)_{t\geq0}$ standard Brownian
motion.
\et

\begin{figure}[ht!]  
\centering%
\subfigure{\includegraphics[width=0.445\textwidth]{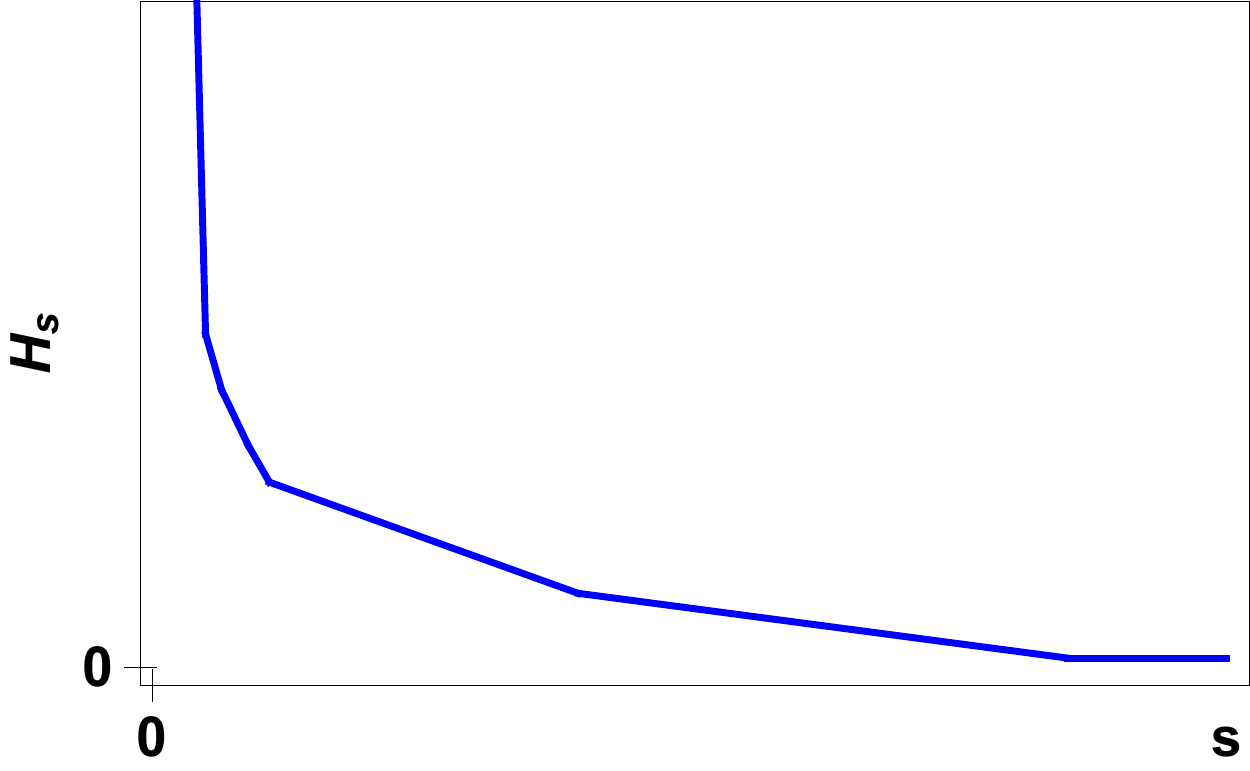}} \:\:\:\:
\subfigure{\includegraphics[width=0.453\textwidth]{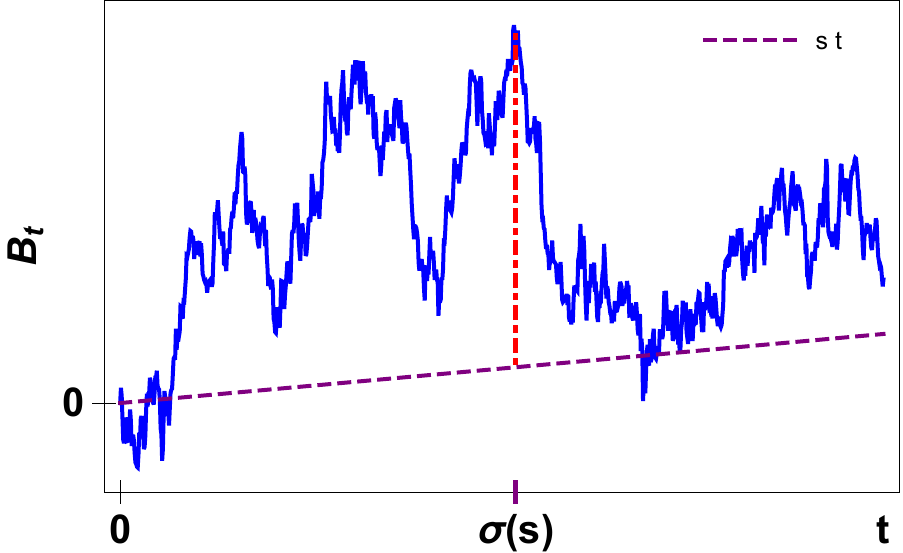}}
\caption{
Simulation showing a typical trajectory of $(H_s)_{s\geq 0}$ (left
panel). Letting $\sig(s):=\inf\{t\geq 0:B_t-st=H_s\}$ be the position where
$B_t-st$ assumes its maximum, we observe that
$H_s=B_{\sig(s)}-\sig(s)s$ and hence $-\dif{s}H_s=\sig(s)$. 
}
\label{fig:H}
\end{figure}

It is known that the function $(H_s)_{s>0}$ of Theorem~\ref{T:main} is
piecewise linear. Let $H'_s:=\dif{s}H_s$ denote the left derivative of $H_s$
and set $\tau(a):=-H'_{1/a}$ $(a>0)$ with $\tau(0):=0$. In \cite{Gro83}, it
is proved that $(\tau(a))_{a\geq 0}$ is a jump process with independent
non-stationary increments. Moreover, its number of jumps in an interval
$(a,b)$ with $0<a<b<\infty$ is Poisson distributed with mean $\log(b/a)$. Note
that this implies that for $0<a<b<\infty$, the number of points in $(a,b)$
where the derivative of $H_s$ makes a jump is also Poisson distributed with
mean $\log(b/a)$. See Figure~\ref{fig:H} for a typical trajectory of
$(H_s)_{s>0}$ and an explanation of the fact that $-H'_s=\inf\{t\geq
0:B_t-st=H_s\}$ $(s>0)$.

\begin{figure}[ht!]  
\centering%
\subfigure{\includegraphics[width=0.446\textwidth]{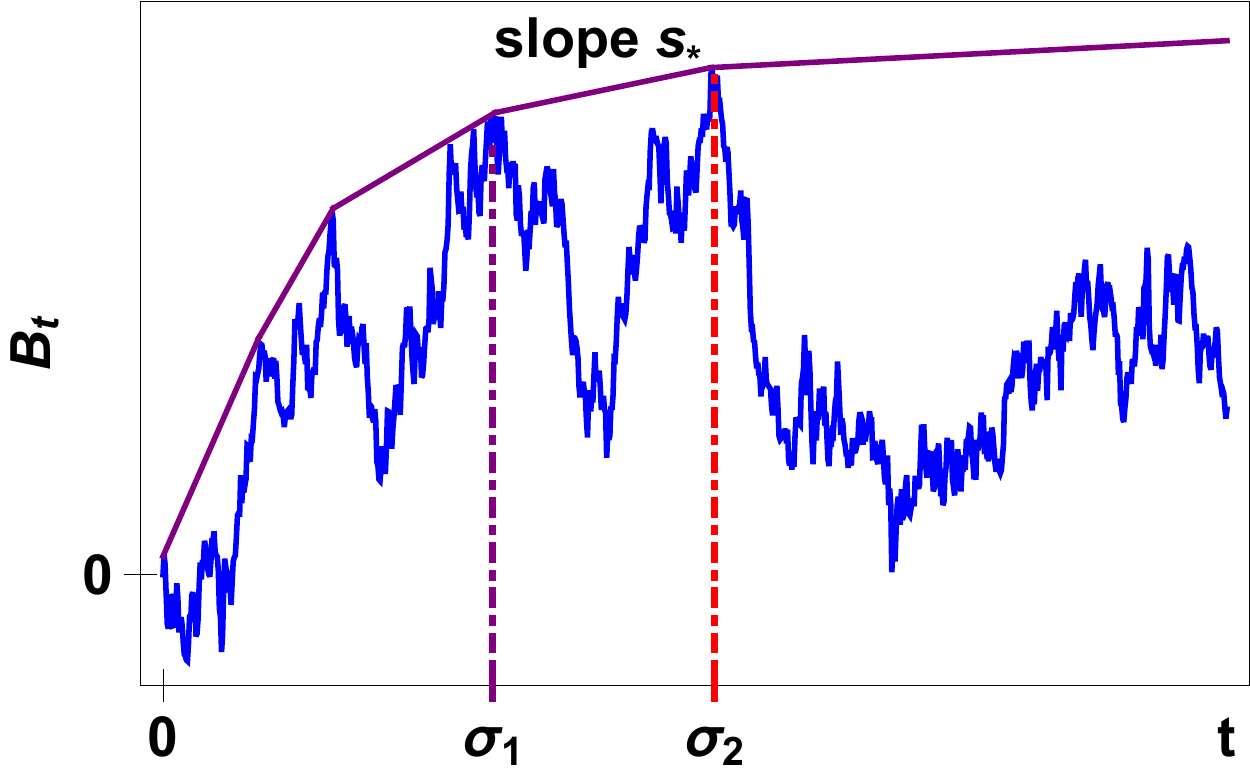}} \:\:\:\:
\subfigure{\includegraphics[width=0.46\textwidth]{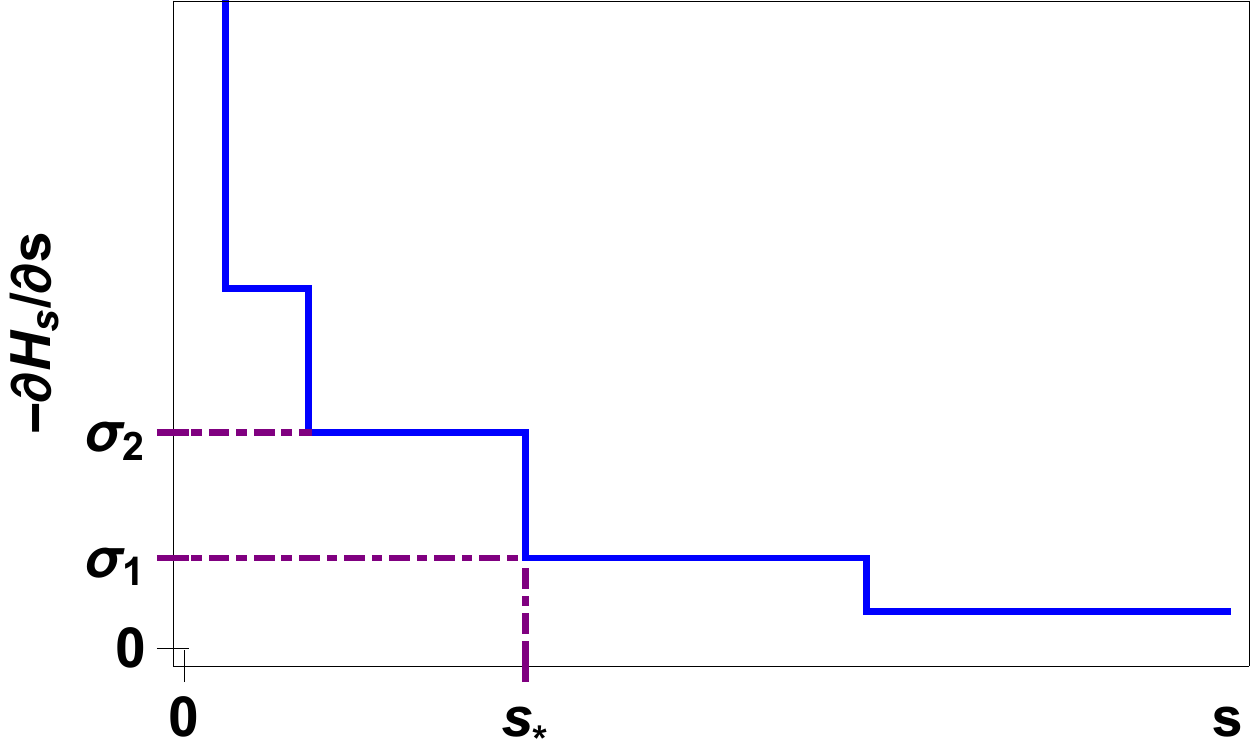}}
\caption{
Relation of $H_s$ to the convex hull of Brownian motion (on the left). The
slope $s_\star$ of a line segment of the convex hull corresponds to a value of
$s$ where the position $\sig(s)$ of the maximum of $B_t-st$ makes a jump,
which corresponds to a jump of $-\dif{s}H_s$ (see Figure~\ref{fig:H}).
}
\label{fig:H2}
\end{figure}

In our model, the quantity $-\dif{s}H_s$ describes the local Poisson density
of tasks with priority close to $s$ and hence is proportional to the time that
has passed since the last time that the highest-priori task $Y_1$ in the queue
had a value below $s$. A jump of $-\dif{s}H_s$ at some priority $s_\star$
reflects the fact that at some quite distant moment in the past, $Y_1$, coming
from the right, reached a local minimum at $s_\star$ before moving up
again. The relation between jumps of $-\dif{s}H_s$ and the convex hull of
Brownian motion is explained in Figure~\ref{fig:H2}.

\subsection{Discussion and overview}

The inbox model of Gabrielli and Caldarelli \cite{CG09} that is our object of
study in the present paper exhibits self-organized criticality. Indeed, the
model organizes itself in such a way that critical behavior associated with
the transition between recurrence and transience can be observed due to
incoming tasks with a priority close to the critical point $\la_{\rm
  c}=1$. This expresses itself in a power law for serving times as demonstated
in \cite{CG09} and also in our main result Theorem~\ref{T:main} that shows
that in equilibrium, the number of tasks with priority above $-\la$ is of
order $(\la_{\rm c}-\la)^{-1}$.

Gabrielli and Caldarelli's model is reminiscent of the well-known Bak Sneppen
model \cite{BS93}, which is one of the best-known models exhibiting
self-organized criticality, although this is only been fully rigorously
established for a simplified version of the model \cite{MS12}. Other similar
models are a ``one-sided canyon model'' introduced in \cite{Swa15} and a model
for traders placing limit buy and sell limit orders at a stock market, first
introduced by Stigler \cite{Sti64} and, in a more general form, by Luckock
\cite{Luc03}. All these models are based on a variation of the rule ``kill the
highest (or lowest) particle'', and all these models exhibit a steady state
where particles cluster near a critical point. We also mention
one-dimensional spatial branching processes where at each branching event, the
lowest particles are killed to keep the number of particles constant. These
models also exhibit self-organized criticality. They have been introduced with
a biological interpretation in \cite{Bru06} and have since also been studied
in the mathematical literature in e.g.\ \cite{Mai13}.

It is interesting to speculate to what degree our main result
Theorem~\ref{T:main} describes universal behavior in this class of models. We
expect the statement to be true for the steady states of the models in
\cite{Swa15} and \cite{Luc03}, and possibly (in a somewhat modified form)
for the model in \cite{MS12}. Proving this will be considerably more difficult
than for the present model, however. A great simplifying property of the inbox
model is that the number of tasks $N^\la(t)$ above a fixed priority is a
Markov process, and the same is true for the restriction of $Y(t)$ to
$[-\la,0]$. For the model in \cite{Swa15}, this first property fails but the
second is still true; for the model in \cite{Luc03}, both fail.

As an open problem for the inbox model, we mention the following. In our main
result Theorem~\ref{T:main}, we rescale space and the weight of items, but we
look only at one fixed time. Is it possible to rescale also time and obtain a
Markov process taking values in the space $\Hi$ that has the law on the
right-hand side of (\ref{toH}) as its invariant law?

The rest of the paper is devoted to proofs. After some initial observations
and definitions in Subsection~\ref{S:lower}, we prove our main result
(Theorem~\ref{T:main}) in Subsection~\ref{S:diffu}. The proof depends on some
lemmas that are proved in Subsection~\ref{S:lemmas}. The paper concludes with
the proof of Proposition~\ref{P:ergod} in Subsection~\ref{S:ergod}.

\section{Proofs}

\subsection{The lower invariant process}\label{S:lower}

Extend the Poisson point sets $\Pi_{\rm in}$ and $\Pi_{\rm out}$ to negative times,
i.e., let $\Pi_{\rm in}$ be a Poisson point set on $(-\infty,0]\times\R$ with
intensity one and $\Pi_{\rm out}$ a Poisson process on $\R$ with intensity
1. For each starting time $s\in\R$ and locally finite subset
$y\sub(-\infty,0]$, we set
\be\label{Ybdef}
\Yb_{s,t}(y):=Y(t)\qquad(t\geq s)
\ee
where $(Y(t))_{t\geq s}$ is the inbox process started at time $s$ in the
initial state $Y(s)=y$ and defined in terms of $\Pi_{\rm in}$ and $\Pi_{\rm
  out}$ as in Section~\ref{S:descrip}. Then $(\Yb_{s,t})_{s\leq t}$ is a
stochastic flow, i.e., a collection of random maps such that $\Yb_{s,s}$ is
the identity map and $\Yb_{t,u}\circ\Yb_{s,t}=\Yb_{s,u}$, almost surely for
all $s\leq t\leq u$. The next lemma says that these maps are monotone with
respect to set inclusion.

\bl[Monotonicity]\label{L:monot}
Almost surely, $y\sub\ti y$ implies $\Yb_{s,t}(y)\sub\Yb_{s,t}(\ti y)$ for all
$s\leq t$ and locally finite subsets $y,\ti y\sub(-\infty,0]$.
\el
\bpro
Let $Y(t):=\Yb_{s,t}(y)$ and $\ti Y(t):=\Yb_{s,t}(\ti y)$. If $t$ is a time
when an incoming task of priority $\la$ arrives, then this task is added both
to $Y$ and $\ti Y$, so $Y(t-)\sub\ti Y(t-)$ implies $Y(t)\sub\ti Y(t)$. If $t$
is a time when a task is executed, then the task with the highest priority (if
there is one) is removed from both to $Y$ and $\ti Y$. If $Y(t-)\sub\ti
Y(t-)$, then either the highest priority element $\ti Y_1(t-)$ is not an
element of $Y(t-)$, or $\ti Y_1(t-)=Y_1(t-)$, so also in this case the
inclusion is preserved.
\epro

The construction of the stationary process $(\un Y(t))_{t\in\R}$ from the next
lemma is similar to the construction of the lower invariant law of a
monotone interacting particle system (see \cite[Thm~III.2.3]{Lig85}).

\bl[Lower invariant process]\label{L:lowinv}
Almost surely, for each $t\in\R$ there exists a random countable subset $\un
Y(t)\sub(-\infty,0]$ such that
\be\label{lowinv}
\Yb_{s,t}(\emptyset)\up\un Y(t)\quad\mbox{as }s\down-\infty.
\ee
\el
\bpro
For each $s\leq s'$, one has
$\Yb_{s',s'}(\emptyset)=\emptyset\sub\Yb_{s,s'}(\emptyset)$. Using the
stochastic flow property and Lemma~\ref{L:monot}, we see that
$\Yb_{s',t}(\emptyset)\sub\Yb_{s,t}(\emptyset)$ for all $s\leq s'\leq t$.
It follows that the left-hand side of (\ref{lowinv}) increases to a limit as
$s\down-\infty$.
\epro

In line with earlier notation, we denote
\be\label{unYN}
\un Y^\la(t):=\un Y(t)\cap[-\la,0]
\quand
\un N^\la(t):=\big|\un Y^\la(t)\big|
\qquad(t\in\R,\ \la\geq 0).
\ee
For $0\leq\de\leq 1$, we will derive a formula for $\un N^{1-\de}(0)$ that is
reminiscent of the definition of the process $H_s$ from Theorem~\ref{T:main}.
As a first step, we prove the following simple lemma.

\bl[Reflected random walk]
For\label{L:Erefl} $0\leq\de\leq 1$ and $u\leq 0$, let
\be\label{Ede}
E^\de(u):=\big|\Pi_{\rm out}\cap[u,0]\big|
-\big|\Pi_{\rm in}\cap\big([-1+\de,0]\times[u,0]\big)\big|
\ee
denote the number of times that a task is executed in the time
interval $[u,0]$ minus the number of tasks with priority $\geq-(1-\de)$ that
arrive in the time interval $[u,0]$. Then
\be\label{Ereflect}
\big|\Yb_{s,u}(\emptyset)\cap[-1+\de,0]\big|
=E^\de(u)-\inf_{s\leq t\leq u}E^\de(t)\qquad(s\leq u\leq 0).
\ee
\el

\bpro
Clearly, the left- and right-hand sides of (\ref{Ereflect}) are both zero if
$s=u$. Increasing $u$ for fixed $s$, we observe that both sides of
(\ref{Ereflect}) increase by one if a task arrives with priority in
$[-(1-\de),1]$. At times of $\Pi_{\rm out}$, either both sides of
(\ref{Ereflect}) are zero and remain zero (due to the fact that both
$E^\de(u)$ and its running infimum decrease by one), or both sides of
(\ref{Ereflect}) are nonzero and decrease by one.
\epro

Setting $u=0$ in (\ref{Ereflect}) and letting $s\down-\infty$, using the fact
that $E^\de(0)=0$, it follows that
\be\label{Esinf}
\un N^{1-\de}(0)=-\inf_{t\leq 0}E^\de(t).
\ee 
In view of what follows, it will be convenient to write the random walk
$E^\de$ as the sum of a driftless random walk (which will converge to Brownian
motion) and a term that contains the drift (which will converge to a linear
function). To this aim, we define, for each $t\geq 0$,
\be\label{Fdef}
F(t):=\big|\Pi_{\rm in}\cap\big([-1,0]\times[-t,0]\big)\big|
-\big|\Pi_{\rm out}\cap[-t,0]\big|
\qquad(t\geq 0).
\ee
In words, $F(t)$ is the number of tasks with priority $\geq-1$ that arrive in
the time interval $[-t,0]$ minus the number of times that a task is executed
in the time interval $[-t,0]$. Next, for $0\leq\de\leq 1$ and $t\geq 0$, we
define
\be\label{Gdef}
G^\de(t):=\big|\Pi_{\rm in}\cap\big([-1,-1+\de]\times[-t,0]\big)\big|
\qquad(t\geq 0),
\ee
which is the number of tasks with priority in $[-1,-(1-\de)]$ that
arrive in the time interval $[-t,0]$.

\bp[Supremum formula]\label{P:supr}
Almost surely, for all $0\leq\de\leq 1$,
\be\label{supr}
\un N^{1-\de}(0)=\sup_{t\geq 0}\big(F(t)-G^\de(t)\big).
\ee
\ep
\bpro
Since $F(t)-G^\de(t)=-E^\de(-t)$ $(t\geq 0)$, formula
(\ref{supr}) is just a rewrite of (\ref{Esinf}).
\epro

\subsection{The diffusive scaling limit}\label{S:diffu}

In this section, we prove our main result Theorem~\ref{T:main}.
In view of (\ref{Heps}), we are interested in $\eps\un N^{1-2\eps s}(0)$, which
by Proposition~\ref{P:supr} is given by
\bc\label{Neps}
\dis\eps\un N^{1-2\eps s}(0)
&=&\sup_{t\geq 0}\big(\eps F(t)-\eps G^{2\eps s}(t)\big)\\[5pt]
&=&\sup_{t\geq 0}\big(\eps F(\ha\eps^{-2}t)-\eps G^{2\eps s}(\ha\eps^{-2}t)\big)
\qquad(\eps>0,\ 0\leq s\leq(2\eps)^{-1}).
\ec
For each $\eps>0$, we define rescaled functions
$F^{(\eps)}:\half\to\R$ and $G^{(\eps)}:\half^2\to\R$ by
\be\left.\ba{r@{\,}c@{\,}l}\label{FGscale}
\dis F^{(\eps)}(t)&:=&\dis\eps F(\ha\eps^{-2}t),\\[5pt]
\dis G^{(\eps)}(s,t)&:=&\dis\eps G^{\,1\wedge 2\eps s}(\ha\eps^{-2}t)
\ea\right\}
\quad(s,t\geq 0).
\ee
We will show that as $\eps\down 0$, the function $F^{(\eps)}$ approximates
Brownian motion and $G^{(\eps)}(s,t)$ approximates $st$. We need this
convergence to be locally uniform in $s$ and $t$. The easiest way to formulate
this is to use \emph{coupling}, i.e., we replace $(F^{(\eps)},G^{(\eps)})$ by
random variables defined on a different underlying probability space, but with
the same distribution as the old ones, so that the convergence is almost sure.

\bl[Convergence of coupled processes]\label{L:coup}
For each $\eps_n\down 0$, it is possible to couple the random variables
$(F^{(\eps_n)},G^{(\eps_n)})$ with $n\geq 0$ in such a way, that almost surely
\be\ba{rrl}\label{coup}
{\rm(i)}\quad&\dis\sup_{t\in[0,T]}\big|F^{(\eps_n)}(t)-B_t\big|\asto{n}0
&\dis\quad\forall T<\infty,\\[5pt]
{\rm(ii)}\quad&\dis
\sup_{(s,t)\in[0,S]\times[0,T]}\big|G^{(\eps_n)}(s,t)-st\big|\asto{n}0
&\dis\quad\forall S,T<\infty,
\ec
where $(B_t)_{t\geq 0}$ is a standard Brownian motion.
\el

The second ingredient in the proof of Theorem~\ref{T:main} is the following
estimate, which guarantees that the supremum over $t\geq 0$ in (\ref{Neps})
and the limit $\eps\down 0$ can be interchanged.

\bl[Uniform upper estimate]\label{L:upper}
For each $s>0$, one has
\be\label{upper}
\lim_{T\to\infty}\sup_{\eps\in(0,1]}
\P\big[F^{(\eps)}(t)-G^{(\eps)}(s,t)\geq 0\mbox{ for some }t\geq T\big]=0.
\ee
\el

The final ingredient for the proof of Theorem~\ref{T:main} is a convergence
criterion for the topology on $\Hi$ if the limit function is
continuous.

\bl[Continuous limit]\label{L:criterion}
Let $\Hi$ be the space of functions defined in Subsection~\ref{S:Hi}.
Let $h^n,h\in\Hi$ and assume that $h$ is continuous.
Then the following statements are equivalent.
\begin{enumerate}
\item $\dis h^n\to h$ in the topology on $\Hi$.
\item $\dis\sup_{s\in[s_0,\infty)}|h^n_s-h_s|\asto{n}0$ for all $s_0>0$.
\item $\dis h^n_s\asto{n}h_s$ for all $s>0$.
\end{enumerate}
\el

We first show how Lemmas~\ref{L:coup}--\ref{L:criterion} imply
Theorem~\ref{T:main}, and in the next subsection then prove the lemmas.\med

\bpro[of Theorem~\ref{T:main}]
Lemma~\ref{L:lowinv} proves the convergence in (\ref{ergod}) in the special
case that $Y(0)=\emptyset$, where $Y^1(\infty)$ is equal in distribution to
$\un Y^1(0)$. Pending the proof of Proposition~\ref{P:ergod}, we will
prove Theorem~\ref{T:main} with the definition (\ref{Heps}) replaced by
$H^\eps_s:=\eps\un N^{1-2\eps s}(0)$. Then (\ref{Neps}) tells us that
\be
H^\eps_s=\sup_{t\geq 0}\big(F^{(\eps)}(t)-G^{(\eps)}(s,t)\big),
\ee
where $F^{(\eps)}$ and $G^{(\eps)}$ are defined in (\ref{FGscale}).

We observe that for any $f,g:[0,T]\to\R$, one has
\be
\big|\sup_{t\in[0,T]}f(t)-\sup_{t\in[0,T]}g(t)\big|\leq\sup_{t\in[0,T]}|f(t)-g(t)|,
\ee
i.e., the map that assigns to a function on $[0,T]$ its supremum is continuous
with respect to the supremum norm. In view of this, let us write
\be\left.\ba{r@{\,}c@{\,}l}
\dis H^\eps_{T,s}&:=&\dis\sup_{t\in[0,T]}\big(F^{(\eps)}(t)-G^{(\eps)}(s,t)\big),\\[5pt]
\dis H_{T,s}&:=&\dis\sup_{t\in[0,T]}\big(B_t-st\big).
\ea\right\}\quad(T<\infty,\ s\geq 0,\ \eps>0).
\ee
Then, for any $\eps_n\down 0$, coupling the random variables
$(F^{(\eps_n)},G^{(\eps_n)})$ as in Lemma~\ref{L:coup}, we observe that
almost surely,
\be\label{Hcutoff}
\sup_{s\in[0,S]}\big|H^{\eps_n}_{T,s}-H_{T,s}\big|\asto{n}0
\qquad(S,T<\infty).
\ee
Fix $s_0>0$. Then
\be
\P\big[B_t-st\geq 0\mbox{ for some }t\geq T,\ s\geq s_0\big]
\leq\P\big[B_t-s_0t\geq 0\mbox{ for some }t\geq T\big]\asto{T}0,
\ee
and hence
\be
\P\big[H_{T,s}\neq H_s\mbox{ for some }s\geq s_0\big]\asto{T}0.
\ee
Similarly, using the fact that $G^{(\eps)}(s,t)$ is nondecreasing in $s$,
we obtain from Lemma~\ref{L:upper} that
\be
\sup_{\eps\in(0,1]}\P\big[H^\eps_{T,s}\neq H^\eps_s
\mbox{ for some }s\geq s_0\big]\asto{T}0.
\ee
Combining this with (\ref{Hcutoff}), we see that
\be
\P\big[\sup_{s\in[s_0,S]}\big|H^{\eps_n}_s-H_s\big|\geq\de\big]
\asto{n}0\qquad(0<s_0<S<\infty,\ \de>0).
\ee
This implies that
\be
\E\big[1\wedge\sup_{s\in[s_0,S]}\big|H^{\eps_n}_s-H_s\big|\big]
\asto{n}0\qquad(0<s_0<S<\infty)
\ee
and hence also
\be
\E\big[\sum_{N=2}^\infty2^{-N}\big(1\wedge\sup_{s\in[1/N,N]}
\big|H^{\eps_n}_s-H_s\big|\big)\big]\asto{n}0.
\ee
It follows that there is a subsequence $\eps_{n(m)}$ such that the expression
in the expectation converges to zero a.s.\ (compare the proof of
\cite[Lemma~3.2]{Kal97}), which by Lemma~\ref{L:criterion} implies that
\be
H^{\eps_{n(m)}}\asto{n}H\quad{\rm a.s.},
\ee
where $\to$ denotes convergence in the topology on $\Hi$. Letting $\Rightarrow$
denote weak convergence of probability laws on $\Hi$ with respect to this
topology, it follows that every sequence $\eps_n\down 0$ contains a subsequence
$\eps_{n(m)}$ such that
\be
\P\big[(H^{\eps_{n(m)}}_s)_{s>0}\in\,\cdot\,\big]
\Asto{m}
\P\big[(H_s)_{s>0}\in\,\cdot\,\big],
\ee
proving (\ref{toH}).
\epro

\subsection{Proof of the lemmas}\label{S:lemmas}

In this section, we complete the proof of Theorem~\ref{T:main} by proving
Lemmas~\ref{L:coup}--\ref{L:criterion}.\med

\bpro[of Lemma~\ref{L:coup}]
Since $F(t)$ is a continuous-time random walk on $\Z$ that jumps one step up
or down with rate one each, the existence of a coupling such that
(\ref{coup})~(i) holds follows by standard arguments, but to also get 
(\ref{coup})~(ii) we have to work a bit.

We observe from (\ref{FGscale}) and (\ref{Gdef}) that
\be
G^{(\eps)}(s,t)
=\eps\big|\Pi_{\rm in}\cap[-1,0\wedge(1-2\eps s)]\times[-\ha\eps^{-2}t,0]\big|.
\ee
Consider the map $\psi_\eps:\R^2\to\R^2$ defined as
\be
\psi_\eps(p,t):=\big((2\eps)^{-1}(1+p),-2\eps^2t\big),
\ee
and let $\Xi^\eps$ denote the random sum of delta measures
\be
\Xi^\eps:=\sum_{(p,t)\in\Pi_{\rm in}}\de_{\psi_\eps(p,t)}.
\ee
Then $\Xi^\eps$ is a Poisson point process on $(-\infty,(2\eps)^{-1}]\times\R$
with intensity $\eps^{-1}$, and
\be\label{GXi}
G^{(\eps)}(s,t)=\eps\int 1_{[0,s]\times[0,t]}\di\Xi^\eps,
\ee
where $1_A$ denotes the indicator function of a set $A$. We can couple the
processes $\Xi^\eps$ for different values of $\eps$ in a monotone way, i.e.,
such that processes with intensities $\eps_1^{-1}\leq\eps_2^{-2}$ satisfy
$\Xi^{\eps_1}\leq\Xi^{\eps_2}$ a.s. (Note that also the domain
$(-\infty,(2\eps)^{-1}]\times\R$ is a monotone function of $\eps^{-1}$.)
Setting up such a coupling in the obvious way, we will have that for each
bounded measurable $A\sub\R^2$, the process $t\mapsto\Xi^{1/t}(A)$ is (for
large enough $t$) a standard Poisson process, i.e., a Markov process on $\N$
that jumps one step up with rate one and never jumps down. The strong law of
large numbers then implies that for this sort of coupling,
\be\label{Xicoup}
\eps\Xi^\eps\Ato{\eps}{0}\ell\quad{\rm a.s.},
\ee
where $\ell$ denotes the Lebesgue measure on $\R^2$ and $\Rightarrow$ denotes
vague convergence of locally finite measures on $\R^2$.

We equip the space of all cadlag (right-continuous with left limits) functions
$f:\half\to\R$ with the Skorohod topology (see \cite[Appendix~A2]{Kal97}) and
we equip the space of all locally finite measures on $\R^2$ with the topology
of vague convergence. Standard results show that with respect to the Skorohod
topology, $F^{(\eps)}=(F^{(\eps)}(t))_{t\geq 0}$ converges weakly in law to
standard Brownian motion as $\eps\down 0$, while (\ref{Xicoup}) shows that
with respect to the topology of vague convergence, $\eps\Xi^\eps$ converges
weakly in law to the deterministic limit $\ell$. We can now apply the Skorohod
representation theorem \cite[Thm~3.30]{Kal97} to conclude that for each
$\eps_n\down 0$, the random variables $(F^{(\eps_n)},\eps_n\Xi^{\eps_n})$ can be
coupled such that a.s.,
\be
F^{(\eps_n)}\asto{n}B
\quand
\eps_n\Xi^{\eps_n}\Asto{n}\ell,
\ee
where $B=(B_t)_{t\geq 0}$ is a standard Brownian motion, $\to$ denotes
convergence w.r.t.\ the Skorohod topology and $\Rightarrow$ denotes vague
convergence of locally finite measures on $\R^2$.

Since $B$ has continuous sample paths, the a.s.\ convergence of $F^{(\eps_n)}$
to $B$ in the Skorohod topology is equivalent to locally uniform convergence,
which gives us (\ref{coup})~(i). Since the indicator functions
$1_{[0,s]\times[0,t]}$ are a.e.\ continuous with respect to Lebesgue measure,
by (\ref{GXi}), the a.s.\ vague convergence of $\eps_n\Xi^{\eps_n}$ to $\ell$
implies the a.s.\ pointwise convergence
\be
G^{(\eps_n)}(s,t)\asto{n}st\qquad(s,t\geq 0).
\ee
Since $G^{(\eps_n)}$ is a.s.\ monotone in both $s$ and $t$,
Lemma~\ref{L:monmon} below allows us to conclude that
this convergence must in fact be uniform on rectangles of the form
$[0,S]\times[0,T]$ as in (\ref{coup})~(ii).
\epro

\bl[Convergence of monotone functions]\label{L:monmon}
Let $C:=I_1\times\cdots\times I_d$ be a hypercube in $[-\infty,\infty]^d$,
where $I_i=[I_i^-,I_i^+]$ is a compact interval for each
$i=1,\ldots,d$. Denote elements of $C$ as $x=(x_1,\ldots,x_d)$ and equip $R$
with the partial order $x\leq y$ iff $x_i\leq y_i$ for all $i=1,\ldots,d$. Let
$g_n,g$ be real functions such that $\lim_{n\to\infty}g_n(x)=g(x)$ for all $x\in C$.
Assume that each $g_n$ is monotone in the sense that $x\leq y$ implies
$g_n(x)\leq g_n(y)$. Then
\be\label{monmon}
\sup_{x\in C}\big|g_n(x)-g(x)\big|\asto{n}0.
\ee
\el
\bpro
We first prove the statement for $C\sub\R^d$. Assume that (\ref{monmon}) does
not hold. Then we can find $\eps>0$ and $x(n)\in C$ such that
$|g_n(x(n))-g(x(n))|\geq\eps$ for all $n$. By the compactness of $C$, going to
a subsequence if necessary, we can assume that $x(n)\to x$ for some $x\in
C$. Let
\be
K_\de(x):=C\cap\{y\in\R^d:|x_i-y_i|\leq\de\ \forall i=1,\ldots,d\}.
\ee
By the continuity of $g$, we can find $\de>0$ such that
$|g(y)-g(x)|\leq\ha\eps$ for all $y\in K_\de(x)$.
Since $K_\de(x)$ is the intersection of two hypercubes, it is itself a
hypercube, i.e., there exist $x^-$ and $x^+$ such that
\be
K_\de(x)=\{y\in\R^d:x^-\leq y\leq x^+\}.
\ee
Since $x(n)\to x$, we have $x(n)\in K_\de(x)$ for all $n$ large enough, and
the monotonicity of $g_n$ then implies that
\be
g_n(x^-)\leq g_n(x(n))\leq g_n(x^+).
\ee
Letting $n\to\infty$, using the pointwise convergence
of $g_n$ to $g$, we see that
\be\ba{l}
\dis g(x)-\ha\eps\leq g(x^-)\leq\liminf_{n\to\infty}g_n(x(n))\\[5pt]
\dis\quad\leq\limsup_{n\to\infty}g_n(x(n))\leq g(x^+)\leq g(x^+)+\ha\eps,
\ec
which contradicts the assumption that $|g_n(x(n))-g(x(n))|\geq\eps$ for
all $n$.

This concludes the proof when $C\sub\R^d$. The more general case
$C\sub[-\infty,\infty]^d$ is not really more general since
$[-\infty,\infty]^d$ is isomorphic to $[0,1]^d$, both in the sense of topology
and in the sense of the partial order $\leq$.
\epro

\bpro[of Lemma~\ref{L:upper}]
Since $G^{(\eps)}(s,t)$ is a.s.\ nondecreasing as a function of $s$, it
suffices to prove the statement for $s$ sufficiently small; in particular, we
can assume without loss of generality that $s\in(0,\ha]$. Then, for each
$\eps\in(0,1]$, by (\ref{FGscale}),
\be
F^{(\eps)}(t)-G^{(\eps)}(s,t)
=\eps F(\ha\eps^{-2}t)-\eps G^{2\eps s}(\ha\eps^{-2}t)\qquad(t\geq 0),
\ee
so we can rewrite (\ref{upper}) as
\be
\lim_{T\to\infty}\sup_{\eps\in(0,1]}
\P\big[F(t)-G^{2\eps s}(t)\geq 0\mbox{ for some }t\geq\ha\eps^{-2}T\big]=0.
\ee
Using notation as in (\ref{Ede}), this says that
\be
\lim_{T\to\infty}\sup_{\eps\in(0,1]}
\P\big[E^{2\eps s}(t)\leq 0\mbox{ for some }t\geq\ha\eps^{-2}T\big]=0.
\ee
Setting $\de:=2\eps s$, this says that
\be
\lim_{T\to\infty}\sup_{\de\in(0,2s]}
\P\big[E^\de(t)\leq 0\mbox{ for some }t\geq2s^2\de^{-2}T\big]=0.
\ee
Since this should hold for any $s\in(0,\ha]$, it suffices to show that
\be
\lim_{T\to\infty}\sup_{\de\in(0,1]}
\P\big[E^\de(t)\leq 0\mbox{ for some }t\geq\de^{-2}T\big]=0,
\ee
where $(E^\de(t))_{t\geq 0}$ is a continuous-time random walk on $\Z$
that starts in $E^\de(0)=0$ and jumps up by one with rate $1$ and down by one
with rate $1-\de$.

Using the well-known fact that
\be
\P\big[E^\de(t)\leq 0\mbox{ for some }t\geq 0\,\big|\,E^\de(0)=x\big]
=1\wedge(1-\de)^x
\ee
and the Markov property, it suffices to show that
\be\label{upper2}
\lim_{T\to\infty}\sup_{\de\in(0,1]}
\E\big[1\wedge(1-\de)^{E^\de(\de^{-2}T)}\big]=0.
\ee
Since $E^\de(t)$ is the difference of two independent Poisson
distributed random variables with mean $t$ and $(1-\de)t$, respectively,
\be
\E\big[E^\de(\de^{-2}T)\big]=\de^{-1}T
\quand
\var\big(E^\de(\de^{-2}T)\big)=(2-\de)\de^{-2}T.
\ee
Estimating $2-\de\leq 2$, Chebyshev's inequality gives
\be
\P\big[|E^\de(\de^{-2}T)-\de^{-1}T|\geq r\de^{-1}\sqrt{2T}\big]\leq r^{-2}
\qquad(r\geq 1).
\ee
It follows that
\be
\E\big[1\wedge(1-\de)^{E^\de(\de^{-2}T)}\big]
\leq(1-\de)^{\de^{-1}(T-r\sqrt{2T})}+r^{-2}\qquad(r\geq 1).
\ee
Since $\log(1-\de)\leq-\de$ on $(0,1]$, we have
\be
(1-\de)^{\de^{-1}}=\ex{\de^{-1}\log(1-\de)}\leq e^{-1}\qquad\big(\de\in(0,1]\big),
\ee
and hence
\be
\E\big[1\wedge(1-\de)^{E^\de(\de^{-2}T)}\big]
\leq\ex{-(T-r\sqrt{2T})}+r^{-2}\qquad(r\geq 1).
\ee
It follows that
\be
\limsup_{T\to\infty}\sup_{\de\in(0,1]}
\E\big[1\wedge(1-\de)^{E^\de(\de^{-2}T)}\big]\leq r^{-2}\qquad(r\geq 1).
\ee
Since $r$ is arbitrary, this proves (\ref{upper2}).
\epro

\bpro[of Lemma~\ref{L:criterion}]
Let $\mu[h]$ denote the locally finite measure on $(0,\infty]$ that $h\in\Hi$
is the distribution function of. Then $\mu[h^n]$ converges vaguely to $\mu[h]$
if and only if $h^n_s\to h_s$ for each $s$ that is a point of continuity of
the limiting function $h$. In particular, if $h$ is continuous, then $h^n\to
h$ in the topology on $\Hi$ if and only if $h^n_s\to h_s$ for all
$s\in(0,\infty)$, proving the equivalence of (i) and (iii). The implication
(ii)$\volgt$(iii) is trivial while (iii)$\volgt$(ii) follows from
Lemma~\ref{L:monmon}. 
\epro

\subsection{Ergodicity}\label{S:ergod}

In this section, we prove Proposition~\ref{P:ergod}. Recall
the definitions of $\un Y^\la(t)$ and $\un N^\la(t)$ from (\ref{unYN}).
Extending our definition in (\ref{Ybdef}), for $\la\geq 0$, we define
random maps $(\Yb^\la_{s,t})_{s\leq t}$ by
\be
\Yb^\la_{s,t}(y):=Y(t)\cap[-\la,0]\qquad(t\geq 0)
\ee
where $y$ is a finite subset of $[-\la,0]$ and $(Y(t))_{t\geq s}$ is the inbox
process started at time $s$ in an initial state with $Y(s)\cap[-\la,0]=y$ and
defined in terms of $\Pi_{\rm in}$ and $\Pi_{\rm out}$ as in
Section~\ref{S:descrip}. Note that by consistency, $Y(t)\cap[-\la,0]$ is a
function of $Y(s)\cap[-\la,0]$ only. 

\bl[Geometric distribution]\label{L:geom}
One has
\be
\P\big[\un N^\la(0)=n\big]=\dis(1-\la)\la^n
\qquad(0\leq\la<1,\ n\geq 0).
\ee
\el
\bpro
We have $\un N^\la(0)=|\un Y^\la(0)|$, where $\un Y^\la(0)$ is the a.s.\ limit
of $\Yb^\la_{-t,0}(\emptyset)$ as $t\to\infty$. Now
$|\Yb^\la_{-t,s}(\emptyset)|$, as a function of $s$, is a Markov process that
jumps as in (\ref{driftwalk}). It is easy to check that for $\la<1$, the
geometric distribution with parameter $\la$ satisfies the detailed balance
conditions and hence is an invariant law; this also proves positive
recurrence. By irredicubility, this is the long-time limit
law started from any initial state, and hence the law of $|\un Y^\la(0)|$.
\epro

\bl[Successful coupling]\label{L:ergod}
For each $0\leq\la<1$ and finite set $y\sub[-\la,0]$, one has
\be\label{ergcoup}
\lim_{t\to\infty}\P\big[\Yb^\la_{-t,0}(y)=\un Y^\la(0)\big]=1.
\ee
\el
\bpro
Since $\un Y^\la(0)$ is a.s.\ finite by Lemma~\ref{L:geom}, and is the a.s.\ limit
of $\Yb^\la_{-t,0}(\emptyset)$ as $t\to\infty$ by its definition in
(\ref{lowinv}), in order to prove (\ref{ergcoup}), it suffices to show that
\be
\lim_{t\to\infty}\P\big[\Yb^\la_{-t,0}(y)=\Yb^\la_{-t,0}(\emptyset)\big]=1.
\ee
By monotonicity (Lemma~\ref{L:monot}),
$\Yb^\la_{-t,s}(y)\supset\Yb^\la_{-t,s}(\emptyset)$ for all $s\in[-t,0]$.
In particular, if $\Yb^\la_{-t,s}(y)=\emptyset$ for some $s\in[-t,0]$,
then $\Yb^\la_{-t,s}(\emptyset)=\emptyset$ for the same $s$ and the two
processes are equal at each later time, so
\be\label{sempt}
\P\big[\Yb^\la_{-t,0}(y)=\Yb^\la_{-t,0}(\emptyset)\big]
\geq\P\big[\Yb^\la_{-t,0}(y)=\emptyset\mbox{ for some }s\in[-t,0]\big].
\ee
Since $|\Yb^\la_{-t,s}(y)|$, as a function of $s$, is a Markov process that
jumps as in (\ref{driftwalk}), which is recurrent for $\la\leq 1$,
the right-hand side of (\ref{sempt}) tends to one as $t\to\infty$.
\epro

\bpro[of Proposition~\ref{P:ergod}]
Defining $Y^1(\infty):=\un Y^1(0)$, the convergence in total variation
distance in (\ref{ergod}) is an immediate consequence of the coupling in
Lemma~\ref{L:ergod}, while $N^\la(\infty)$ is geometrically distributed by
Lemma~\ref{L:geom}. Since $|Y^\la(\infty)|$ is geometrically distributed with
parameter $\la$, we see that $|Y^\la(\infty)|$ tends to infinity in
probability as $\la\up 1$. Since $Y^\la(\infty)\up Y^1(\infty)$, this implies
that $Y^1(\infty)$ is a.s.\ an infinite set.
\epro

\section*{Acknowledgments}

Work sponsored by GA\v{C}R grant P201/12/2613.



\bibliographystyle{plain}

\end{document}